\documentclass[technote]{IEEEtran}

\usepackage{amsmath,graphicx}
\newtheorem{theorem}{\bf Theorem}
\newtheorem{proposition}{\bf Proposition}
\newtheorem{lemma}{\bf Lemma}
\newtheorem{corollary}{\bf Corollary}

\begin{document}
\title{On Log-concavity of the Generalized Marcum Q Function}

\author{Yaming~Yu, {\it Member, IEEE}
\thanks{Yaming Yu is with the Department of Statistics, University of California, Irvine, CA, 92697-1250, USA
(e-mail: yamingy@uci.edu).  This work is supported in part by a start-up fund from the Bren School of Information and 
Computer Sciences at UC Irvine.}
}

\maketitle

\begin{abstract}
It is shown that, if $\nu\geq 1/2$ then the generalized Marcum Q function $Q_\nu(a, b)$ is log-concave in $b\in [0, \infty)$.  This proves a conjecture of Sun, Baricz and Zhou (2010).  We also point out relevant results in the statistics literature. 
\end{abstract}

\begin{keywords}
increasing failure rate; log-concavity; modified Bessel function; noncentral chi square. 
\end{keywords}

\section{Introduction}
The generalized Marcum Q function \cite{M60} has important applications in radar detection and communications over fading channels and has received much attention; see, e.g., \cite{BS09, H92, KMK09}, \cite{LK06}-\cite{S89} and \cite{SA03}-\cite{SBZ}.  It is defined as 
\begin{equation}
\label{def1}
Q_\nu(a, b)=\int_b^\infty \frac{t^\nu}{a^{\nu-1}} \exp\left(-\frac{t^2+a^2}{2}\right) I_{\nu-1}(at)\, {\rm d}t
\end{equation}
where $\nu>0,\ a, b\geq 0$ and $I_\nu$ denotes the modified Bessel function of the first kind of order $\nu$ defined by the series \cite{AS72} (9.6.10)
$$I_\nu(t) = \sum_{k=0}^\infty \frac{(t/2)^{2k+\nu}}{k!\Gamma(\nu+k+1)}.$$
($Q_\nu(0, b)$ is defined by taking $a\downarrow 0$.)  Recently, Sun, Baricz and Zhou \cite{SBZ} have studied the monotonicity, log-concavity, and tight bounds of $Q_\nu(a,b)$ in great detail.  We are concerned with log-concavity, which has intrinsic interest, and can help establish useful bounds; see \cite{SBZ} and the references therein for the large literature in information theory and communications on numerical calculations of $Q_\nu(a,b)$. 

This note resolves some of the conjectures made by \cite{SBZ}.  We also point out relevant literature in statistics on both
theoretical properties and numerical computation of $Q_\nu(a,b)$.  Our Theorem~\ref{thm1} proves Conjecture~1 of \cite{SBZ}. 
\begin{theorem}
\label{thm1}
The function $Q_\nu(a, b)$ is log-concave in $b\in [0, \infty)$ for all $a\geq 0$ if and only if $\nu\geq 1/2$. 
\end{theorem}

A sufficient condition for log-concavity of an integral like (\ref{def1}) is that the integrand is log-concave in $t$.  Proposition~\ref{prop1} and Theorem~\ref{hard} take this approach. 

\begin{proposition} 
\label{prop1}
The integrand in (\ref{def1}) is log-concave in $t\in (0, \infty)$ for all $\nu\geq 1/2$ if and only if $0\leq a\leq 1$. 
\end{proposition}

\begin{theorem}
\label{hard}
The integrand in (\ref{def1}) is log-concave in $t\in (0,\infty)$ for all $a\geq 0$ if and only if $\nu\geq \nu_0$ where $\nu_0\approx 0.78449776$ is the unique solution of the equation 
$$\frac{I_\nu(\sqrt{5-2\nu})}{I_{\nu-1}(\sqrt{5-2\nu})}=\frac{3-2\nu}{\sqrt{5-2\nu}}$$
in the interval $\nu\in (1/2, 3/2)$. 
\end{theorem}

Note the difference between Proposition~\ref{prop1} and Theorem~\ref{hard}: the former gives a criterion for log-concavity in $t$ for all $\nu\geq 1/2$ whereas the latter gives one for all $a\geq 0$.  From Proposition~\ref{prop1} and Theorem~\ref{hard} we obtain Corollary~\ref{coro1}, which confirms part of Conjecture~2 of \cite{SBZ}. 
\begin{corollary}
\label{coro1}
The function $1-Q_\nu(a,b)$ is log-concave in $b\in [0,\infty)$, if either (i) $\nu\geq 1/2$ and $0\leq a\leq 1$, or (ii) $\nu\geq \nu_0$ as in Theorem~\ref{hard}. 
\end{corollary}

The case of $Q_1(a, b)$ (Marcum's original Q function) is especially interesting.  If $\nu=1$ then the integrand in (\ref{def1}) is the probability density function (PDF) of a Rice distribution, $Q_1(a,b)$ being the corresponding tail probability, or survival function.  Therefore Theorem~\ref{hard} yields 
\begin{corollary}
The probability density function, cumulative distribution function (CDF), and survival function of a Rice distribution are all log-concave. 
\end{corollary}

In general, let $X$ be a noncentral $\chi^2$ random variable with $2\nu$ degrees of freedom and noncentrality parameter $a^2$.  Then
$$Q_\nu(a,b)=\Pr\left(\sqrt{X}>b\right).$$
Equivalently, $1-Q_\nu(\sqrt{a}, \sqrt{b})$ is the CDF of a noncentral $\chi^2$ random variable with $2\nu$ degrees of freedom and noncentrality parameter $a$.  The noncentral $\chi^2$ distribution plays an important role in statistical hypothesis testing and has been extensively studied.  We mention \cite{D92, KB96} on numerical computation and \cite{FR97, JKB, S79} on theoretical properties.  Its CDF, and hence $Q_\nu(a, b)$, can be routinely calculated (e.g., using ${\tt pchisq()}$ in the R package). 

Concerning theoretical properties, Finner and Roters \cite{FR97} (see also \cite{DGS84}) have obtained the following results  using tools from total positivity \cite{K68}. 
\begin{theorem}[\cite{FR97}, Theorems 3.4, 3.9; Remark 3.6]
\label{thm2}
The function $1-Q_\nu(\sqrt{a}, \sqrt{b})$ is log-concave 
\begin{itemize}
\item
in $b\in[0, \infty)$ for $\nu>0,\, a\geq 0$;
\item
in $\nu>0$ for $a,\, b\geq 0$;
\item
in $a\geq 0$ for $\nu>0,\, b\geq 0$.
\end{itemize}
The function $Q_\nu(\sqrt{a}, \sqrt{b})$ is log-concave
\begin{itemize}
\item
in $b\in [0, \infty)$ for $\nu\geq 1,\, a\geq 0$;
\item
in $\nu\in [1/2, \infty)$ for $a,\, b\geq 0$;
\item
in $a\geq 0$ for $\nu>0,\, b\geq 0$. 
\end{itemize}
\end{theorem}
Theorem~\ref{thm2} and Corollary~\ref{coro1} cover several results of \cite{SBZ}, including part of their Conjectures~2 and 3 (see also \cite{SB08}).  The parts of these conjectures that remain open are 
\begin{itemize}
\item
$1-Q_\nu(a, b)$ is log-concave in $b\in [0,\infty)$ for $\nu\in [1/2, \nu_0)$ and $a> 1$; 
\item
$Q_\nu(a, b)$ is log-concave in $\nu\in (0, 1/2]$ for $a,\, b\geq 0$. 
\end{itemize}

In Section~II we prove Theorems~\ref{thm1}, \ref{hard} and Proposition~\ref{prop1}.  The proof of Theorem~\ref{thm1} uses a general technique which may be helpful in related problems.  The proof of Theorem~\ref{hard} relies partly on numerical verification as theoretical analysis appears quite cumbersome. 

\section{Proof of Main Results} 
The following observation, which is of independent interest, is key to our proof of Theorem~\ref{thm1}. 

\begin{lemma}
\label{lem1}
Let $f(t)$ be a probability density function on $\mathbf{R}\equiv (-\infty, \infty)$. Assume (i) $f(t)$ is unimodal, i.e., there exists $t_0\in \mathbf{R}$ such that $f(t)$ increases on $(-\infty, t_0]$ and decreases on $[t_0, \infty)$; (ii) $f(t_0-)\leq f(t_0+)$; (iii) $f(t)$ is log-concave in the declining phase $t\in (t_0, \infty)$.  Then the survival function $\bar{F}(b)\equiv \int_b^\infty f(t)\, {\rm d}t$ is log-concave in $b\in \mathbf{R}$. 
\end{lemma}
\begin{IEEEproof}
Assumption (iii) implies that $\bar{F}(b)$ is log-concave in $b\in [t_0,\infty)$.  Because $f(t)$ increases on $(-\infty, t_0]$, we know $\bar{F}(b)$ is concave and hence log-concave on $(-\infty, t_0]$.  By Assumption (ii) we have 
$$\bar{F}'(t_0-) =-f(t_0-)\geq -f(t_0+)=\bar{F}'(t_0+).$$ 
Hence $\bar{F}(b)$ is log-concave in $b\in \mathbf{R}$ overall. 
\end{IEEEproof}

{\bf Remark 1.} A distribution whose survival function is log-concave is said to have an increasing failure rate (IFR) \cite{BP75}.  Distributions with IFR form an important class in reliability and survival analysis.  Lemma~\ref{lem1} provides a simple sufficient condition for IFR distributions. 

Henceforth let $f(t)$ be the integrand in (\ref{def1}) for $t>0$.  Equivalently, $f(t)$ is the density function of a noncentral $\chi$ random variable with $2\nu$ degrees of freedom.  Define
\begin{equation}
\label{r}
r_\nu(t)=\frac{I_\nu(t)}{I_{\nu-1}(t)}.
\end{equation}
We use $r'_\nu(t)$ to denote the derivative with respect to $t$. 
\begin{lemma}
\label{lem2}
If $\nu\geq 1/2$ then $f'(t)/(tf(t))$ decreases in $t\in (0, \infty)$.
\end{lemma}
\begin{IEEEproof}
Let us assume $\nu>1/2$ and $a>0$.  The boundary cases follow by taking limits.  Direct calculation yields 
\begin{align}
\nonumber
\frac{f'(t)}{tf(t)} &= \frac{\nu}{t^2} - 1 + \frac{aI_{\nu-1}'(at)}{tI_{\nu-1}(at)}\\
\label{deri}
&= \frac{2\nu-1}{t^2} - 1 + \frac{a r_{\nu}(at)}{t}
\end{align}
where (\ref{deri}) uses (\ref{r}) and the formula \cite{AS72} (9.6.26)
\begin{equation}
\label{recur}
I_{\nu-1}'(t) = I_\nu(t) + \frac{\nu-1}{t} I_{\nu-1}(t).
\end{equation}
Since $(2\nu-1)/t^2$ decreases in $t$, we only need to show that $r_\nu(t)/t$ decreases in $t$.  We may use the integral formula of \cite{AS72} (9.6.18) and obtain 
$$\frac{r_\nu(t)}{t} = \frac{\int_0^1 (1-s^2) g(s, t)\, {\rm d} s} {(2\nu-1)\int_0^1 g(s, t)\, {\rm d} s}$$
where 
$$g(s, t)=(1-s^2)^{\nu-3/2} \cosh(ts).$$  
As can be easily verified, if $0<t_1<t_2$ then $g(s, t_2)/g(s, t_1)$ increases in $s\in (0, 1)$.  That is, $g(s, t)$ is TP$_2$ \cite{K68}.  Since $1-s^2$ decreases in $s\in (0,1)$, by Proposition 3.1 in Chapter~1 of \cite{K68}, the ratio $\int_0^1 (1-s^2) g(s, t)\, {\rm d} s/ \int_0^1 g(s, t)\, {\rm d} s$ decreases in $t\in (0, \infty)$, as required. 
\end{IEEEproof}
\begin{IEEEproof}[Proof of Theorem~\ref{thm1}]
Let us assume $\nu>1/2$ and show log-concavity.  By Lemma~\ref{lem2}, either (i) $f'(t)< 0$ for all $t\in (0, \infty)$ or (ii) there exists some $t_0\in (0, \infty)$ such that $f'(t)\geq 0$ when $t< t_0$ and $f'(t)\leq 0$ when $t> t_0$.  (Since $\int_0^\infty f(t)\, {\rm d} t = Q_\nu(a, 0)=1$, it cannot happen that $f'(t)> 0$ for all $t\in (0, \infty)$.)  In either case $f(t)$ satisfies Assumptions (i) and (ii) of Lemma~\ref{lem1} ($f(t)\equiv 0$ for $t\leq 0$).  Let us consider Case (ii); the same argument applies to Case (i).  For $t\in (t_0, \infty)$ we have $f'(t)\leq 0$, and hence  
\begin{align*}
\frac{1}{t} \frac{{\rm d}^2}{{\rm d}t^2} \log f(t) &\leq \frac{1}{t} \frac{{\rm d}^2}{{\rm d}t^2} \log f(t) - \frac{f'(t)}{t^2 f(t)}\\
&= \frac{\rm d}{{\rm d}t} \left(\frac{1}{t} \frac{\rm d}{{\rm d} t} 
\log f(t)\right) \leq 0
\end{align*}
where the last step holds by Lemma~\ref{lem2}.  Thus $f(t)$ is log-concave in $t\in (t_0, \infty)$ and Assumption (iii) of Lemma~\ref{lem1} is satisfied.  We conclude that $Q_\nu(a, b)=\int_b^\infty f(t)\, {\rm d}t$ is log-concave in $b\in [0, \infty)$.

It remains to show that, if $Q_\nu(a, b)$ is log-concave in $b\in [0, \infty)$ for all $a\geq 0$, then we must have $\nu\geq 1/2$.  Let us consider $a=0$.  We have
$$Q_\nu(0, b)=1-\frac{1}{2^\nu\Gamma(\nu)}\int_0^{b^2} t^{\nu-1} e^{-t/2}\, {\rm d}t.$$
As $b\downarrow 0$, it is easy to see that $\log Q_\nu(0, b)$ behaves like 
$$\log\left(1- C b^{2\nu} + o(b^{2\nu})\right) = -C b^{2\nu} +o(b^{2\nu})$$
with $C=2^{-\nu}/\Gamma(\nu+1)$.  Hence, if $\nu< 1/2$ then $Q_\nu(0, b)$ is no longer log-concave for $b$ near zero.  It follows that the $1/2$ in Theorem~\ref{thm1} is the best possible. 
\end{IEEEproof}

\begin{IEEEproof}[Proof of Proposition~\ref{prop1}]
Using (\ref{deri}) we get 
\begin{equation}
\label{l2}
\frac{{\rm d}^2}{{\rm d}t^2}\log f(t) =-\frac{2\nu-1}{t^2} -1 + a^2 r'_\nu(at).
\end{equation}
However, 
\begin{align}
\nonumber
r'_\nu(t) &= \frac{I'_\nu(t)}{I_{\nu-1}(t)}-\frac{I_\nu(t)I_{\nu-1}'(t)}{I_{\nu-1}^2(t)}\\
\label{rprime}
 & = 1-\frac{2\nu-1}{t} r_\nu(t) - r^2_\nu(t)
\end{align}
where (\ref{rprime}) holds by applying (\ref{r}), (\ref{recur}) and the recursion \cite{AS72} (9.6.26) 
$$I_{\nu+1} (t) =I_{\nu-1}(t) - \frac{2\nu}{t} I_\nu(t).$$
If $\nu\geq 1/2$ and $0<a\leq 1$ then $r'_\nu(at)\leq 1$ by (\ref{rprime}), and we have  
$$\frac{{\rm d}^2}{ {\rm d}t^2} \log f(t)\leq a^2 - 1 \leq 0.$$
Hence $f(t)$ is log-concave in $t\in (0, \infty)$. 

To show the converse, suppose $f(t)$ is log-concave in $t$ for all $\nu\geq 1/2$.  Consider $\nu=1/2$.  As $t\downarrow 0$ we have $r_\nu(t)\to 0$, and ${\rm d}^2 \log f(t)/{\rm d}t^2 \to a^2-1$.  Hence we must have $a\leq 1$. 
\end{IEEEproof}

{\bf Remark 2.} For $\nu\geq 1/2$, the function $f(t)$ is log-concave in its declining phase, as shown in the proof of Theorem~\ref{thm1}.  If $a\in [0,1]$ in addition, then Proposition~\ref{prop1} shows that $f(t)$ is log-concave in all $t\in (0, \infty)$.  For $a>1$ and $\nu\geq 1/2$, however, numerical evidence suggests that $f(t)$ may not be log-concave in its rising phase.  Hence a version of Lemma~\ref{lem1} cannot be applied to $1-Q_\nu(a, b)$.  Log-concavity of $1-Q_\nu(a,b)$ in $b$ appears to be a difficult problem. 

Let us establish two lemmas before proving Theorem~\ref{hard}.
\begin{lemma}
\label{lem3}
The function $f(t)$ is log-concave in $t\in (0, \infty)$ for all $a\geq 0$ if and only if the function 
\begin{equation}
\label{ht}
h_\nu(t)=1-\frac{2\nu-1}{t^2} -\frac{2\nu-1}{t}r_\nu(t) - r^2_\nu(t)
\end{equation}
is nonpositive for $t\in (0, \infty)$. 
\end{lemma} 
\begin{IEEEproof}
By (\ref{rprime}) we get 
\begin{equation}
\label{ht2}
h_\nu(t)=r'_\nu(t) - \frac{2\nu-1}{t^2}.
\end{equation}
If $h_\nu(t)\leq 0$ then by (\ref{l2}) we have 
\begin{equation*}
\frac{{\rm d}^2}{{\rm d} t^2}\log f(t) = a^2 h_\nu(at) -1 < 0.
\end{equation*}
Conversely, if $f(t)$ is log-concave in $t\in (0, \infty)$ for all $a\geq 0$, then holding $at$ constant while letting $a\to \infty$ yields $h_\nu(s)\leq 0$ for each $s\in (0, \infty)$. 
\end{IEEEproof}

\begin{lemma}
\label{lem4}
The function 
\begin{equation*}
r_\nu(\sqrt{5-2\nu}) - \frac{3-2\nu}{\sqrt{5-2\nu}}
\end{equation*}
strictly increases in $\nu\in [1/2, 3/2]$ and has a zero at $\nu_0\approx 0.78449776$. 
\end{lemma}
\begin{IEEEproof}
Although this only involves a one-variable function over a small interval, it is verified by numerical calculations, as theoretical analysis becomes complicated.  The value of $\nu_0$ is computed by a fixed point algorithm. 
\end{IEEEproof}

\begin{IEEEproof}[Proof of Theorem~\ref{hard}]
Define $h_\nu(t)$ as in (\ref{ht}) and $\nu_0$ as in Lemma~\ref{lem4}.  We examine the intervals $(0, 1/2]$, $(1/2, \nu_0)$ and $[\nu_0, \infty)$ for $\nu$ in turn.  If $0<\nu\leq 1/2$ then letting $t\downarrow 0$ we have $r_\nu(t)\to 0$ and $h_\nu(t)>0$ for small $t$.  By Lemma~\ref{lem3}, $f(t)$ is not log-concave for all $a\geq 0$. 

Let us assume $\nu>1/2$.  Differentiating (\ref{ht}) with respect to $t$ and applying (\ref{ht2}) we get 
\begin{align}
\label{hprime1}
h'_\nu(t) = -\frac{2\nu-1}{t^2}l_\nu(t) - \left(\frac{2\nu-1}{t} + 2r_\nu(t)\right) h_\nu(t)
\end{align}
where 
\begin{equation}
\label{lt}
l_\nu(t)=r_\nu(t) -\frac{3-2\nu}{t}.
\end{equation}
For $\nu>1/2$ we know $r_\nu(t)$ increases from $0$ to $1$ as $t$ increases from $0$ to $\infty$ (see \cite{A74}).  Hence, if $1/2<\nu< 3/2$, then $l_\nu(t)$ strictly increases and $l_\nu(t)=0$ has a unique solution, say at $t_1\in (0, \infty)$.  If $1/2<\nu<\nu_0$, then by Lemma~\ref{lem4}, $l_\nu(\sqrt{5-2\nu})< 0$, and hence $t_1> \sqrt{5-2\nu}$.  In view of (\ref{ht}) and (\ref{lt}) we have
\begin{align}
\label{eqn1}
h_\nu(t_1) &= 1-\frac{2\nu-1}{t_1^2} -\frac{2\nu-1}{t_1}\left(\frac{3-2\nu}{t_1}\right) - \frac{(3-2\nu)^2}{t_1^2}\\
\label{eqn2}
& = 1-\frac{5-2\nu}{t_1^2} > 0.
\end{align}
By Lemma~\ref{lem3}, $f(t)$ is no longer log-concave for all $a\geq 0$. 

Suppose $\nu> \nu_0$.  We have $h_\nu(t)\to -\infty$ as $t\downarrow 0$ and $h_\nu(t)\to 0$ as $t\to\infty$.  If $h_\nu(t)$ does become positive, then there exists a finite $t_0>0$ such that $h_\nu(t_0)=0$ and $h'_\nu(t_0)\geq 0$ (at least one sign change should be from $-$ to $+$).  We get $l_\nu(t_0)\leq 0$ from (\ref{hprime1}).  If $\nu\geq 3/2$ then (\ref{lt}) yields  $l_\nu(t_0)\geq r_\nu(t_0) > 0$, a contradiction.  Hence $h_\nu(t)\leq 0$ for all $t\in (0, \infty)$ if $\nu\geq 3/2$. 

Suppose $\nu_0< \nu<3/2$.  If $l_\nu(t_0)=h'_\nu(t_0)=0$ then we deduce $t_0=\sqrt{5-2\nu}$ from (\ref{ht}) and (\ref{hprime1}) by a calculation similar to (\ref{eqn1})--(\ref{eqn2}).  But $l_\nu(\sqrt{5-2\nu})=0$ contradicts Lemma~\ref{lem4}.  Hence we may assume $h'_\nu(t_0)>0$ and $l_\nu(t_0)<0$.  By Lemma~\ref{lem4} we have $l_\nu(\sqrt{5-2\nu})> 0$.  Because $l_\nu(t)$ is strictly increasing, and $t_1$ is the solution of $l_\nu(t)=0$, we obtain $t_0< t_1<\sqrt{5-2\nu}$.  The calculation (\ref{eqn1})--(\ref{eqn2}) now yields $h_\nu(t_1)<0$.  Because $h_\nu(t_0)=0,\ h'_\nu(t_0)> 0$ there exists $t_*\in (t_0, t_1)$ such that $h_\nu(t_*)=0$ and $h'_\nu(t_*)\leq 0$.  By (\ref{hprime1}), we get $l_\nu(t_*)\geq 0$, which contradicts the strict monotonicity of $l_\nu(t)$ as $l_\nu(t_1)=0$.  It follows that $h_\nu(t)\leq 0,\ t\in (0,\infty),$ and $f(t)$ is log-concave.  Taking the limit we extend this log-concavity to $\nu=\nu_0$. 
\end{IEEEproof}



\end{document}